\documentclass[10pt]{article}
\usepackage{amsfonts}
\usepackage{amssymb}
\newcommand{\fidemo}{\vrule height4pt width3pt depth2pt}
\newcommand{\rees}[1]{\mbox{$\mathcal{R}(#1)$}}
\newcommand{\reesw}[2]{\mbox{$\mathcal{R}(#1;#2)$}}
\newcommand{\agr}[1]{\mbox{$\mathcal{G}(#1)$}}
\newcommand{\agrw}[2]{\mbox{$\mathcal{G}(#1;#2)$}}

\newtheorem{rteorema}{Theorem}
\newtheorem{teorema}{Theorem}[section]

\newtheorem{definicio}[teorema]{Definition}
\newtheorem{proposicio}[teorema]{Proposition}
\newtheorem{lema}[teorema]{Lemma}
\newtheorem{corollari}[teorema]{Corollary}
\newtheorem{observacio}[teorema]{Remark}
\newtheorem{exemple}[teorema]{Example}
\newtheorem{exemples}[teorema]{Examples}
\newtheorem{notacions}[teorema]{Notations}

\begin{document}

\title{The strong uniform Artin-Rees property in codimension one}

\author{{\sc Francesc Planas-Vilanova} \\ \\ {\small
Dept. Matem\`atica Aplicada I.  ETSEIB-UPC.  Diagonal 647, 08028
Barcelona.}  
\addtocounter{footnote}{-1}{\ }\thanks{\hspace*{-17pt}e-mail:
planas@ma1.upc.es} 
}

\date{}

\maketitle

\section{Introduction}

The purpose of this paper is to prove the following theorem:

\begin{rteorema}\label{cucu}
Let $A$ be an excellent (in fact $J-2$) ring and let $N\subseteq M$ be
two finitely generated $A$-modules such that ${\rm dim}(M/N)\leq 1$.
Then there exists an integer $s\geq 1$ such that, for all integers
$n\geq s$ and for all ideals $I$ of $A$,
\begin{eqnarray*}
I^{n}M\cap N=I^{n-s}(I^{s}M\cap N)\, .
\end{eqnarray*}
\end{rteorema}

This result is a variation of a theorem of Duncan and O'Carroll
\cite{do}: maximal ideals are replaced for any ideal using the
unavoidable hypothesis ${\rm dim}(M/N)\leq 1$ (as an Example of Wang
shows \cite{wang1}). Moreover it provides a partial positive answer to
the question raised by Huneke in Conjecture 1.3 \cite{huneke1}.

We begin by recalling what is called uniform Artin-Rees
properties. Let $A$ be a noetherian ring, $I$ be an ideal of $A$ and
let $N\subseteq M$ be two finitely generated $A$-modules. The usual
Artin-Rees lemma states that there exists an integer $s\geq 1$,
depending on $N$, $M$ and $I$, such that for all $n\geq s$,
\begin{eqnarray*}
I^{n}M\cap N=I^{n-s}(I^{s}M\cap N)\, .
\end{eqnarray*}
In particular, $I^{n}M\cap N\subseteq I^{n-s}N$. As in \cite{huneke1},
let us say the pair $(N,M)$ has the ({\em strong}) {\em uniform
Artin-Rees property} with respect to a set of ideals $\mathcal{W}$ of
$A$ and with ({\em strong}) {\em uniform number} $s$ ($s$ depending on
$(N,M;\mathcal{W})$) if, for every ideal $I$ of $\mathcal{W}$ and for
all $n\geq s$, ($I^{n}M\cap N=I^{n-s}(I^{s}M\cap N)$) $I^{n}M\cap
N\subseteq I^{n-s}N$. Clearly, if $s$ is a (strong) uniform number for
$(N,M,\mathcal{W})$ and $t\geq s$, then $t$ is also a (strong) uniform
number for $(N,M,\mathcal{W})$. The minimum of all such (strong)
uniform numbers will be denoted by $s=s(N,M,\mathcal{W})$ and call it
``the'' (strong) uniform number for $(N,M,\mathcal{W})$. If
$\mathcal{W}$ is the set of all ideals of $A$, we delete the phrase
``with respect to $\mathcal{W}$'' and simply write $s=s(N,M)$.

Eisenbud and Hochster \cite{eh} ask whether a pair $(N,M)$ has the
uniform Artin-Rees property with respect to the set of maximal ideals
of $A$. O'Carroll \cite{ocarroll1} proves that if $A$ is excellent
then $A$ has the uniform Artin-Rees property with respect to the set
of maximal ideals and Duncan and O'Carroll \cite{do} generalize this
result to the strong uniform Artin-Rees property. Later, O'Carroll
\cite{ocarroll2} shows the strong uniform Artin-Rees property with
respect to the set of principal ideals of a noetherian ring $A$.
Nevertheless, the strong uniform Artin-Rees property cannot hold for
the class of all ideals of $A$. Indeed, Wang \cite{wang1} shows that
if $(A,\mathfrak{m})$ is a 3-dimensional regular local ring,
$\mathfrak{m}=(x,y,z)$, $I_{k}=(x^{k},y^{k},x^{k-1}y+z^{k})$ and
$J=(z)$, then there does not exist an $s\geq 1$ such that, for all
$n\geq s$ and for all $k\geq 1$, $I^{n}_{k}\cap
J=I^{n-s}_{k}(I^{s}_{k}\cap J)$. Remark that ${\rm dim}(A/J)=2$. Thus,
in this sense, Theorem \ref{cucu} is not improvable.

On the other hand Huneke \cite{huneke1} shows the uniform Artin-Rees
property with respect to the class of all ideals of a noetherian ring
$A$ if $A$ is either essentially of finite type over a noetherian
local ring, either a ring of characteristic $p$ and a module finite
over $A^{p}$ or either essentially of finite type over
$\mathbb{Z}$. In the same paper Huneke conjectures that this theorem
remains true for excellent noetherian rings of finite Krull
dimension. Thus Theorem \ref{cucu} gives a partial positive answer to
this conjecture.

Since strong uniform Artin-Rees property is not true in general,
Huneke \cite{huneke2} asks for classes of ideals where strong uniform
Artin-Rees property holds. If $A$ is regular local, $J$ is an ideal of
$A$, does there exist an $s\geq 1$ such that, for all $n\geq s$ and
for all ideal $I$ of $A$ whose image in $A/J$ is generated by a system
of parameters, $I^{n}\cap J=I^{n-s}(I^{s}\cap J)$ ? In fact, Lai
\cite{lai} proves that this property is equivalent to the
Relation-Type Conjecture, stated by Huneke, and proved by Wang
\cite{wang2} for rings with finite local cohomology. The {\em relation
type} of an ideal $I$, ${\rm rt}(I)$, is the largest degree of any
minimal homogeneous system of generators of the ideal defining the
Rees algebra of $I$. The Relation-Type Conjecture asks whether there
is an integer $s\geq 1$ such that, for all parameter ideal $I$ of a
complete local equidimensional noetherian ring $A$, the relation type
of $I$ is ${\rm rt}(I)\leq s$.

In order to prove Theorem \ref{cucu} we generalise this relationship
between the strong uniform Artin-Rees property and the existence of
uniform bounds for the relation type. First we define ${\rm rt}(I;M)$,
the relation type of an ideal $I$ with repect to an $A$-module $M$
(Section \ref{ja}). Then we consider ${\rm grt}(M)={\rm sup}\{ {\rm
rt}(I;M)\mid I\, \mbox{ideal} \; \mbox{of}\ A\} $, the supremum
(possibly infinite) of all relation types of ideals $I$ of $A$ with
respect to $M$, and call it the {\em global relation type} of the
$A$-module $M$. We prove:

\begin{rteorema}
Let $A$ be a commutative ring, $\mathcal{W}$ a set of ideals of $A$,
$I\in \mathcal{W}$ and $N\subseteq M$ two $A$-modules. Let
$s(N,M;\mathcal{W})$ denote the strong uniform number for the pair
$(N,M)$ with respect to $\mathcal{W}$. Then $s(N,M;\{ I\} )\leq {\rm
rt}(I;M/N)\leq {\rm max}({\rm rt}(I;M),s(N,M;\{ I\}))$. In particular,
$s(N,M;\mathcal{W})\leq {\rm sup}\{ {\rm rt}(J;M/N)\mid J\in
\mathcal{W}\}$ and $s(N,M)\leq {\rm grt}(M/N)$.
\end{rteorema}

We thus ask for whether a module has finite global relation type. A
very special is already known: for a commutative (non necessarily
noetherian local) domain $A$, ${\rm grt}(A)=1$ is equivalent to $A$ be
a ring of Pr\"ufer \cite{costa} and, more in general, commutative
rings with ${\rm grt}(A)=1$ are known to be the rings of weak
dimension one or less \cite{planas1}. Thus, for a noetherian local
ring $A$, ${\rm grt}(A)=1$ if and only if $A$ is a discrete valuation
ring or a field.

Our guide here is the following celebrated theorem of Cohen and Sally
\cite{cohen}, \cite{sally}: for a commutative noetherian local ring
$(A,\mathfrak{m},k)$, ${\rm sup}\{ \mu (I)\mid I \mbox{ ideal of } A\}
<\infty $ is equivalent to dimension of $A$ be ${\rm dim}\, A\leq 1$,
where $\mu (I)={\rm dim}_{k}(I/\mathfrak{m}I)$ stands for the minimum
number of generators of $I$. Then, we show the expected analoguos
result by replacing $\mu (I)$ for the relation type ${\rm rt}(I)$ of
$I$. Concretely:

\begin{rteorema}\label{main}
Let $A$ be an excellent (in fact $J-2$) ring. The following conditions
are equivalent:
\begin{itemize}
\item[$(i)$] ${\rm grt}(M)<\infty $ for all finitely generated
$A$-module $M$.
\item[$(ii)$] ${\rm grt}(A)<\infty $.
\item[$(iii)$] There exists an $r\geq 1$ such that ${\rm
rt}(I)\leq r$ for every three-generated ideal $I$ of $A$.
\item[$(iv)$] There exists an $r\geq 1$ such that $(x^{r}y)^{r}\in
(x^{r+1},y^{r+1})(x^{r+1},y^{r+1},x^{r}y)^{r-1}$ for all $x,y\in A$.
\item[$(v)$] ${\rm dim}\, A\leq 1$.
\end{itemize}
\end{rteorema}

The paper is organized as follows. Section \ref{ja} is dedicated to
recall some definitions and properties on the module of effective
relations of a graded algebra. In order to prove Theorem 2, we need to
generalize them from graded algebras to graded modules. Once
introduced all the machinery, we prove Theorem 2 in Section
\ref{rtim}. In Section \ref{1dim} we prove that rings of finite global
relation type have dimension one or less and in Section \ref{0dim} we
show that zero dimensional modules over noetherian rings have finite
global relation type. This is half of the proof in Theorem 3. In
Section \ref{jor}, we first consider the local case and reduce to
Cohen-Macaulay modules. Then we give a new proof, now for modules, of
a well known result for rings (see, for instance, \cite{trung2}): if
$I$ is an $\mathfrak{m}$-primary ideal of a one dimensional
Cohen-Macaulay local ring $A$ and $M$ is a maximal Cohen-Macaulay
$A$-module, then ${\rm rt}(I;M)\leq e(A)$, the relation type of $I$
with respect $M$ is bounded above by the multiplicity of $A$. We
conclude that one dimensional finitely generated modules over
noetherian local rings have finite global relation type.  Section
\ref{final} finishes with all proofs.  Throughout, $A$ denotes a
commutative ring with unity.  All tensor products are over $A$ unless
specified the contrary. Dimension of a ring or module always mean
Krull dimension. One of the main tools in this note is the module of
effective relations of a graded algebra or module. In order to recall
some of their general properties we will often refer to
\cite{planas2}.

\section{Preliminaries}\label{ja}

Let $A$ be a commutative ring. By a {\em standard} $A$-algebra we mean
a commutative graded $A$-algebra $U=\oplus _{n\geq 0}U_{n}$ with
$U_{0}=A$ and such that $U$ is generated as an $A$-algebra by the
elements of $U_{1}$. Put $U_{+}=\oplus _{n>0}U_{n}$ the {\em
irrelevant ideal} of $U$.

If $E=\oplus _{n\geq 0}E_{n}$ is a graded $U$-module and $r\geq 0$ is
an integer, we denote by $F_{r}(E)$ the submodule of $E$ generated by
the elements of degree at most $r$. Put (possibly infinite) $s(E)={\rm
min}\{ r\geq 1\mid E_{n}=0 \mbox{ for all } n\geq r+1\} $. Remark that
we are only interested for $s(E)\geq 1$. Since for all $n\geq 1$,
$(E/U_{+}E)_{n}=E_{n}/U_{1}E_{n-1}$, then for all $r\geq 1$, the
following three conditions are equivalent: $F_{r}(E)=E$;
$s(E/U_{+}E)\leq r$; and $E_{n}=U_{1}E_{n-1}$ for all $n\geq r+1$.

If $f:V\rightarrow U$ is a surjective graded morphism of standard
$A$-algebras, we denote by $E(f)$ the graded $A$-module $E(f)={\rm
ker}f/V_{+}{\rm ker}f=\oplus _{n\geq 1}{\rm ker}f_{n}/V_{1}{\rm
ker}f_{n-1}=\oplus _{n\geq 1}E(f)_{n}$. The following is an elementary
but very useful fact (Lemma 2.1 \cite{planas2}): if $f:V\rightarrow U$
and $g:W\rightarrow V$ are two surjective graded morphisms of standard
$A$-algebras, then there exists a graded exact sequence of $A$-modules
$E(g)\rightarrow E(f\circ g)\buildrel g\over \rightarrow
E(f)\rightarrow 0$. In particular, $s(E(f))\leq s(E(f\circ g))\leq
{\rm max}(s(E(f)),s(E(g)))$. Moreover, if $V$ and $W$ are two
symmetric algebras, then $E(g)_{n}=0$ and $E(f\circ g)_{n}=E(f)_{n}$
for all $n\geq 2$.

Let $U$ be a standard $A$-algebra, let ${\bf S}(U_{1})$ be the
symmetric algebra of $U_{1}$ and let $\alpha :{\bf
S}(U_{1})\rightarrow U$ be the surjective graded morphism of standard
$A$-algebras induced by the identity on $U_{1}$. The {\em module of
effective} $n$-{\em relations} of $U$ is defined to be
$E(U)_{n}=E(\alpha )_{n}={\rm ker}\alpha _{n}/U_{1}{\rm ker}\alpha
_{n-1}$ (for $n=0,1$, $E(U)_{n}=0$). Put $E(U)=\oplus _{n\geq
2}E(U)_{n}=\oplus _{n\geq 2}E(\alpha )_{n}=E(\alpha )= {\rm ker}\alpha
/{\bf S}_{+}(U_{1}){\rm ker}\alpha$. The {\em relation type} of $U$ is
defined to be ${\rm rt}(U)=s(E(U))$, that is, ${\rm rt}(U)$ is the
minimum positive integer $r\geq 1$ such that the effective
$n$-relations are zero for all $n\geq r+1$.

A {\em symmetric presentation} of $U$ is a surjective graded morphism
of standard $A$-algebras $f:V\rightarrow U$, where $V={\bf S}(V_{1})$
is the symmetric $A$-algebra of the $A$-module $V_{1}$ (for instance,
$V_{1}=U_{1}$ and $f_{1}=1$, or $f_{1}:V_{1}\rightarrow U_{1}$ a free
presentation of $U_{1}$). Using Lemma 2.1 in \cite{planas2} one
deduces that $E(U)_{n}=E(f)_{n}$ for all $n\geq 2$ and
$s(E(U))=s(E(f))$. Thus the module of effective $n$-relations and the
relation type of a standard $A$-algebra are independent of the chosen
symmetric presentation.

For an ideal $I$ of $A$, the module of effective $n$-relations and the
relation type of $I$ are defined to be $E(I)_{n}=E(\rees{I})_{n}$ and
${\rm rt}(I)={\rm rt}(\rees{I})$, where $\rees{I}=\oplus _{n\geq
0}I^{n}t^{n}\subset A[t]$ is the {\em Rees algebra of} $I$. If ${\rm
rt}(I)<\infty $ (for instance, if $A$ is noetherian), then ${\rm
rt}(I)={\rm rt}(\agr{I})$, where $\agr{I}=\oplus _{n\geq
0}I^{n}/I^{n+1}$ is the {\em associated graded ring of} $I$
(Proposition 3.3 \cite{planas2}).

\medskip

Let us now extend the classical notion of relation type of an ideal to
the relation type of an ideal with respect to a module. Some of the
results we present here are a straightforward generalization of former
results. We thus will skip some details.

\begin{definicio}{\rm 
Let $U=\oplus _{n\geq 0}U_{n}$ be a standard $A$-algebra and $F=\oplus
_{n\geq 0}F_{n}$ a graded $U$-module. We will say $F$ is a {\em
standard} $U$-module if $F$ is generated as an $U$-module by the
elements of $F_{0}$, that is, $F_{n}=U_{n}F_{0}$ for all $n\geq 0$.  In
particular, $F_{n}=U_{1}F_{n-1}$ for all $n\geq 1$. 
}\end{definicio}

\begin{exemples}{\rm Some of the most interesting standard modules for
our purposes are the following:
\begin{itemize}
\item[$(1)$] If $U$ is a standard $A$-algebra and $M$ is an
$A$-module, then $U\otimes M$ is a standard $U$-module ($M$ in degree
zero). If $U_{1}=A^{\oplus n}$ is a finitely generated free $A$-module
and $U={\bf S}(U_{1})$ is the symmetric algebra of $U_{1}$, then
$U\otimes M=A[T_{1},\ldots ,T_{n}]\otimes M=M[T_{1},\ldots ,T_{n}]$.
\item[$(2)$] $\reesw{I}{M}=\oplus _{n\geq 0}I^{n}M$, the {\em Rees module
of an ideal} $I$ {\em of} $A$ {\em with respect to an} $A$-{\em
module} $M$, is a standard $\rees{I}$-module.
\item[$(3)$] $\agrw{I}{M}=\oplus _{n\geq 0}I^{n}M/I^{n+1}M$, the {\em
associated graded module of an ideal} $I$ {\em of} $A$ {\em with
respect to an} $A$-{\em module} $M$, is a standard $\agr{I}$-module.
\end{itemize}
}\end{exemples}

Let $U=\oplus _{n\geq 0}U_{n}$ be a standard $A$-algebra and $F=\oplus
_{n\geq 0}F_{n}$, $G=\oplus _{n\geq 0}G_{n}$ be two graded
$U$-modules. If $\varphi :G\rightarrow F$ is a surjective
graded morphism of $U$-modules, we denote by $E(\varphi )$ the graded
$A$-module $E(\varphi )={\rm ker}\varphi /U_{+}{\rm ker}\varphi ={\rm
ker}\varphi _{0}\oplus (\oplus _{n\geq 1}{\rm ker}\varphi
_{n}/U_{1}{\rm ker}\varphi _{n-1})=\oplus _{n\geq 0}E(\varphi )_{n}$.
The following is a generalization of Lemma 2.1 in \cite{planas2}:

\begin{lema}\label{gener}
If $\varphi :G\rightarrow F$ and $\psi :H\rightarrow G$ are two
surjective graded morphisms of graded $U$-modules, then there exists a
graded exact sequence $E(\psi )\rightarrow E(\varphi \circ \psi
)\buildrel \psi \over \rightarrow E(\varphi )\rightarrow 0$ of
$A$-modules. In particular, $s(E(\varphi ))\leq s(E(\varphi \circ \psi
))\leq {\rm max}(s(E(\varphi )),s(E(\psi )))$. Moreover, if $H={\bf
S}(P)\otimes Q$ is the tensor product of the symmetric algebra of the
$A$-module $P$ with the $A$-module $Q$, $G={\bf S}(M)\otimes N$ is the
tensor product of the symmetric algebra of the $A$-module $M$ with the
$A$-module $N$ and $\psi =f\otimes h$ where $f:{\bf S}(P)\rightarrow
{\bf S}(M)$ is induced by an epimorphism $f_{1}:P\rightarrow M$ and
$h:Q\rightarrow N$ is also surjective, then $E(\psi )_{n}=0$ and
$E(\varphi \circ \psi )_{n}=E(\varphi )_{n}$ for all $n\geq 2$.
\end{lema}

\noindent {\em Proof}. To deduce the existence of the exact sequence
we proceed as in Lemma 2.1 in \cite{planas2}. For the second
assertion, consider the following commutative diagram of exact rows:

\begin{picture}(330,125)(-5,0)

\put(200,100){\makebox(0,0){\mbox{\footnotesize
${\bf \Lambda}_{2}(U_{1})\otimes {\bf S}_{n-2}(P)\otimes Q$}}}
\put(300,100){\makebox(0,0){\mbox{\footnotesize
${\bf \Lambda}_{2}(U_{1})\otimes {\bf S}_{n-2}(M)\otimes N$}}}
\put(370,100){\makebox(0,0){$0$}}

\put(245,100){\vector(1,0){8}}
\put(347,100){\vector(1,0){13}}

\put(100,60){\makebox(0,0){\mbox{\footnotesize
$U_{1}\otimes {\rm ker}\psi _{n-1}$}}}
\put(200,60){\makebox(0,0){\mbox{\footnotesize
$U_{1}\otimes {\bf S}_{n-1}(P)\otimes Q$}}}
\put(300,60){\makebox(0,0){\mbox{\footnotesize
$U_{1}\otimes {\bf S}_{n-1}(M)\otimes N$}}}
\put(370,60){\makebox(0,0){$0$}}

\put(135,60){\vector(1,0){22}}
\put(240,60){\vector(1,0){20}}
\put(340,60){\vector(1,0){20}}

\put(40,20){\vector(1,0){30}}
\put(135,20){\vector(1,0){30}}
\put(235,20){\vector(1,0){30}}
\put(325,20){\vector(1,0){35}}

\put(30,20){\makebox(0,0){$0$}}
\put(100,20){\makebox(0,0){{\footnotesize ${\rm ker}\psi _{n}$}}}
\put(200,20){\makebox(0,0){{\footnotesize ${\bf S}_{n}(P)\otimes Q$}}}
\put(300,20){\makebox(0,0){{\footnotesize ${\bf S}_{n}(M)\otimes N$}}}
\put(370,20){\makebox(0,0){$0$}}

\put(100,50){\vector(0,-1){20}}
\put(200,50){\vector(0,-1){20}}
\put(200,50){\vector(0,-1){16}}
\put(300,50){\vector(0,-1){20}}
\put(300,50){\vector(0,-1){16}}
\put(170,40){\makebox(0,0){{\footnotesize $\partial _{1,n}^{P}\otimes
1_{Q}$}}} 
\put(330,40){\makebox(0,0){{\footnotesize $\partial _{1,n}^{M}\otimes
1_{N}$}}}

\put(200,90){\vector(0,-1){20}}
\put(300,90){\vector(0,-1){20}}
\put(170,80){\makebox(0,0){{\footnotesize $\partial _{2,n}^{P}\otimes
1_{Q}$}}} 
\put(330,80){\makebox(0,0){{\footnotesize $\partial _{2,n}^{M}\otimes
1_{N}$}}}

\put(250,110){\makebox(0,0){\mbox{\footnotesize $1\otimes 
\psi _{n-2}$}}}
\put(250,70){\makebox(0,0){\mbox{\footnotesize $1\otimes 
\psi _{n-1}$}}}
\put(250,10){\makebox(0,0){\mbox{\footnotesize $\psi _{n}$}}}
\end{picture}

\noindent where $\partial _{2,n}^{P}((x\wedge y)\otimes z)=y\otimes
xz-x\otimes yz$ and $\partial _{1,n}^{P}(x\otimes t)=xt$, $x,y\in
U_{1}$, $z\in {\bf S}_{n-2}(P)$, $t\in {\bf S}_{n-1}(P)$ ($\partial
^{M}$ defined analogously). By Theorem 2.4 in \cite{planas2}, the
right and middle columns are exact sequences for all $n\geq 2$ and, by
the snake lemma, ${\rm ker}(\partial _{1,n}^{P}\otimes
1_{Q})\rightarrow {\rm ker}(\partial _{1,n}^{M}\otimes
1_{N})\rightarrow E(\psi )_{n}\rightarrow 0$ are exact for all $n\geq
2$. Since $1\otimes \psi _{n-2}$ is surjective, $E(\psi )_{n}=0$ for
all $n\geq 2$. \fidemo

\medskip

\begin{definicio}{\rm 
Let $U$ be a standard $A$-algebra and $F$ be a standard $U$-module.
Let ${\bf S}(U_{1})$ be the symmetric algebra of $U_{1}$ and let
$\alpha :{\bf S}(U_{1})\rightarrow U$ be the surjective graded
morphism of standard $A$-algebras induced by the identity on $U_{1}$.
Let $\gamma :{\bf S}(U_{1})\otimes F_{0}\buildrel \alpha \otimes
1\over \rightarrow U\otimes F_{0}\rightarrow F$ be the composition of
$\alpha \otimes 1$ with the structural morphism. Since $F$ is a
standard $U$-module, $\gamma $ is a surjective graded morphism of
graded ${\bf S}(U_{1})$-modules. The {\em module of effective}
$n$-{\em relations} of $F$ is defined to be $E(F)_{n}=E(\gamma )_{n}=
{\rm ker}\gamma _{n}/U_{1}{\rm ker}\gamma _{n-1}$ (for $n=0$,
$E(F)_{n}=0$). Put $E(F)=\oplus _{n\geq 1}E(F)_{n}=\oplus _{n\geq
1}E(\gamma )_{n}=E(\gamma )= {\rm ker}\gamma /{\bf S}_{+}(U_{1}){\rm
ker}\gamma$. The {\em relation type} of $F$ is defined to be ${\rm
rt}(F)=s(E(F))$, that is, ${\rm rt}(F)$ is the minimum positive
integer $r\geq 1$ such that the effective $n$-relations are zero for
all $n\geq r+1$.

A {\em symmetric presentation} of a standard $U$-module $F$ is a
surjective graded morphism of standard $V$-modules $\varphi
:G\rightarrow F$, with $\varphi :G=V\otimes M\buildrel f\otimes h\over
\rightarrow U\otimes F_{0}\rightarrow F$ where $f:V\rightarrow U$ is a
symmetric presentation of the standard $A$-algebra $U$,
$h:M\rightarrow F_{0}$ is an epimorphism of $A$-modules and $U\otimes
F_{0}\rightarrow F$ is the structural morphism. Using Lemma
\ref{gener}, one deduces that $E(F)_{n}=E(\varphi )_{n}$ for all
$n\geq 2$ and $s(E(F))=s(E(\varphi ))$. Thus the module of effective
$n$-relations and the relation type of a standard $U$-module are
independent of the chosen symmetric presentation.

For an ideal $I$ of $A$ and an $A$-module $M$, the module of effective
$n$-relations and the relation type of $I$ with repect to $M$ are
defined to be $E(I;M)_{n}=E(\reesw{I}{M})_{n}$ and ${\rm rt}(I;M)={\rm
rt}(\reesw{I}{M})$.  }\end{definicio}

\begin{observacio}\label{prougranm}
{\rm The following are simple, but useful remarks:
\begin{itemize}
\item[$(1)$] If $U$ is a standard $A$-algebra, then $U$ is a standard
$U$-module. Moreover the modules of effective $n$-relations of $U$ as a
standard $A$-algebra and as a standard $U$-module are equal
$E_{A\mbox{\footnotesize -alg}}(U)_{n}=E_{U\mbox{\footnotesize
-mod}}(U)_{n}$. Thus ${\rm rt}_{A\mbox{\footnotesize -alg}}(U)={\rm
rt}_{U\mbox{\footnotesize -mod}}(U)$. In particular, if $I$ is an
ideal of $A$, then $E(I;A)_{n}=E(I)_{n}$ and ${\rm rt}(I;A)={\rm
rt}(I)$.
\item[$(2)$] If $f:V\rightarrow U$ is a surjective graded morphism of
standard $A$-algebras and $F$ is a standard $U$-module, then $F$ is a
standard $V$-module. Moreover $E_{U\mbox{\footnotesize
-mod}}(F)_{n}=E_{V\mbox{\footnotesize -mod}}(F)_{n}$ for all $n\geq 2$
and ${\rm rt}_{U\mbox{\footnotesize -mod}}(F)={\rm
rt}_{V\mbox{\footnotesize -mod}}(F)$.
\item[$(3)$] If $\varphi :R\rightarrow A$ is a surjective homomorphism
of rings, $U$ is standard $A$-algebra and $F$ is a standard
$U$-module, then $V=R\oplus U_{+}$ is a standard $R$-algebra and $F$
is a standard $V$-module. Moreover $E(U)_{n}=E(V)_{n}$ for all $n\geq
2$ and ${\rm rt}(U)={\rm rt}(V)$. Analogously $E_{U\mbox{\footnotesize
-mod}}(F)_{n}=E_{V\mbox{\footnotesize -mod}}(F)_{n}$ for all $n\geq 2$
and ${\rm rt}_{U\mbox{\footnotesize -mod}}(F)={\rm
rt}_{V\mbox{\footnotesize -mod}}(F)$.
\item[$(4)$] If $\varphi :G\rightarrow F$ is a surjective graded
morphism of standard $U$-modules such that ${\rm ker}\varphi _{n}=0$
for all $n\geq t$, then $E(G)_{n}=E(F)_{n}$ for all $n\geq t+1$ and
${\rm rt}(G)\leq \mbox{max}({\rm rt}(F),t)$. If $t=1$, then ${\rm
rt}(G)={\rm rt}(F)$. For instance, if $I$ and $J$ are two ideals of
$A$, ${\rm rt}(I/I\cap J)={\rm rt}((I+J)/J)$.
\item[$(5)$] Let $F$ be a standard $U$-module, $\underline{x}=\{
x_{i}\} $ a (possibly infinite) set of generators of the $A$-module
$U_{1}$ and $\underline{T}=\{ T_{i}\} $ a set of as many variables
over $A$ as $\underline{x}$ has elements. Take $V_{1}=\oplus
_{i}AT_{i}$, $V={\bf S}(V_{1})=A[\underline{T}]$, $G=V\otimes
F_{0}=F_{0}[\underline{T}]$ and $\varphi :G\rightarrow F$ defined by
$\varphi (\sum y_{i}T_{i})=\sum x_{i}y_{i}$. Clearly $\varphi $ is a
symmetric presentation of $F$. Thus, ${\rm rt}(F)=1$ if and only if
${\rm ker}\varphi $ is generated by linear forms. If
$I=(\underline{x})$ is an ideal of $A$ and $M$ is an $A$-module, then
${\rm rt}(I;M)=1$ if and only if the kernel of the surjective graded
morphism $\varphi :M[\underline{T}]\rightarrow \reesw{I}{M}$, $\varphi
(\sum y_{i}T_{i})=\sum x_{i}y_{i}$, is generated by linear forms. We
say $I$ is an {\em ideal of linear type with respect to} $M$ if ${\rm
rt}(I;M)=1$ (\cite{hsv} pag 106, \cite{trung1} pag 41).
\end{itemize}
}\end{observacio}

\noindent {\em Proof}. $(1)$ follows from definitions. $(2)$ is
consequence of Lemma \ref{gener}. For the proof of $(3)$, consider
$\alpha _{V}:{\bf S}^{R}(U_{1})\rightarrow V$ surjective graded
morphism of standard $R$-algebras, $\alpha _{U}:{\bf
S}^{A}(U_{1})\rightarrow U$ surjective graded morphism of standad
$A$-algebras, $f:V\rightarrow U$ and $g:{\bf S}^{R}(U_{1})\rightarrow
{\bf S}^{A}(U_{1})$ the natural surjective graded morphisms extending
$\varphi $. Since $f\circ \alpha _{V}=\alpha _{U}\circ g$ and $f_{n}$
and $g_{n}$ are isomorphisms for all $n\geq 1$, then
$E(V)_{n}=E(\alpha _{V})_{n}=E(f\circ \alpha )_{n}=E(\alpha
_{U})_{n}=E(U)_{n}$. For the rest of $(3)$ is sufficient to apply the
tensor product $-\otimes F_{0}$. In order to prove $(4)$, let $\psi
:H\rightarrow G$ be a symmetric presentation of $G$. By hypothesis,
${\rm ker}\psi_{n}={\rm ker}(\varphi \circ \psi )_{n}$ for all $n\geq
t$. Hence $E(G)_{n}=E(\psi )_{n}=E(\varphi \circ \psi )_{n}=E(F)_{n}$
for all $n\geq t+1$. If $n\geq t+1,{\rm rt}(F)+1$, then
$E(G)_{n}=E(F)_{n}=0$. Thus ${\rm rt}(G)\leq {\rm max}\left( {\rm
rt}(F), t\right) $.  Take $G=\rees{I/I\cap J}$, $F=\rees{(I+J)/J}$ and
$\varphi :G\rightarrow F$ the natural surjective graded morphism with
$\varphi _{0}:A/I\cap J\rightarrow A/J$ and $\varphi _{n}:(I^{n}+I\cap
J)/(I\cap J)\buildrel \simeq \over \rightarrow (I^{n}+J)/J$
isomorphism for all $n\geq 1$.  Applying consecutively $(1)$, $(4)$,
$(3)$ and $(1)$, $E_{A/I\cap J\mbox{\footnotesize -alg}}(G)_{n}=
E_{G\mbox{\footnotesize -mod}}(G)_{n}= E_{G\mbox{\footnotesize
-mod}}(F)_{n}= E_{A/J\mbox{\footnotesize -alg}}(F)_{n}$. Finally,
$(5)$ follows from definitions. \fidemo

\medskip

Let us now modify Theorem 2.4 in \cite{planas2} to modules:

\begin{proposicio}\label{febrer}
Let $U$ be a standard $A$-algebra and let $F$ be a standard
$U$-module. For each integer $n\geq 2$, there exists a complex of
$A$-modules
\begin{eqnarray*}
{\bf \Lambda}_{2}(U_{1})\otimes F_{n-2}\buildrel \partial _{2,n}\over
\longrightarrow U_{1}\otimes F_{n-1}\buildrel \partial _{1,n}\over
\longrightarrow F_{n}\, ,
\end{eqnarray*}
defined by $\partial _{2,n}((x\wedge y)\otimes z)=y\otimes xz-x\otimes
yz$ and $\partial _{1,n}(x\otimes t)=xt$ and whose homology is
$E(F)_{n}$. 
\end{proposicio}

\noindent {\em Proof}. By Theorem 2.4, there exists ${\bf \Lambda}
_{2}(U_{1})\rightarrow U_{1}\otimes U_{1}\rightarrow U_{2}\rightarrow
0$, a complex of $A$-modules defined by $\partial _{2}(x\wedge
y)=y\otimes x-x\otimes y$ and $\partial _{1}(x\otimes t)=xt$. Applying
the tensor product $-\otimes F_{n-2}$ and considering the structural
morphisms $U_{i}\otimes F_{j}\rightarrow F_{i+j}$ we get the complex.
Let ${\bf S}(U_{1})$ be the symmetric algebra of $U_{1}$ and let
$\alpha :{\bf S}(U_{1})\rightarrow U$ be the surjective graded
morphism of standard $A$-algebras induced by the identity on $U_{1}$.
Let $\gamma :{\bf S}(U_{1})\otimes F_{0}\buildrel \alpha \otimes
1\over \rightarrow U\otimes F_{0}\rightarrow F$ be the composition of
$\alpha \otimes 1$ with the structural morphism. Consider now, for
each $n\geq 2$, the following commutative diagram of exact rows:

\begin{picture}(330,125)(-5,0)

\put(200,100){\makebox(0,0){\mbox{\footnotesize
${\bf \Lambda}_{2}(U_{1})\otimes {\bf S}_{n-2}(U_{1})\otimes F_{0}$}}}
\put(300,100){\makebox(0,0){\mbox{\footnotesize
${\bf \Lambda}_{2}(U_{1})\otimes F_{n-2}$}}}
\put(370,100){\makebox(0,0){$0$}}

\put(250,100){\vector(1,0){15}}
\put(330,100){\vector(1,0){30}}

\put(100,60){\makebox(0,0){\mbox{\footnotesize
$U_{1}\otimes {\rm ker}\gamma _{n-1}$}}}
\put(200,60){\makebox(0,0){\mbox{\footnotesize
$U_{1}\otimes {\bf S}_{n-1}(U_{1})\otimes F_{0}$}}}
\put(300,60){\makebox(0,0){\mbox{\footnotesize
$U_{1}\otimes F_{n-1}$}}}
\put(370,60){\makebox(0,0){$0$}}

\put(135,60){\vector(1,0){22}}
\put(240,60){\vector(1,0){30}}
\put(330,60){\vector(1,0){30}}

\put(40,20){\vector(1,0){30}}
\put(135,20){\vector(1,0){30}}
\put(235,20){\vector(1,0){35}}
\put(325,20){\vector(1,0){35}}

\put(30,20){\makebox(0,0){$0$}}
\put(100,20){\makebox(0,0){{\footnotesize ${\rm ker}\gamma _{n}$}}}
\put(200,20){\makebox(0,0){{\footnotesize ${\bf S}_{n}(U_{1})
\otimes F_{0}$}}}
\put(300,20){\makebox(0,0){{\footnotesize $F_{n}$}}}
\put(370,20){\makebox(0,0){$0$}}

\put(100,50){\vector(0,-1){20}}
\put(200,50){\vector(0,-1){20}}
\put(200,50){\vector(0,-1){16}}
\put(300,50){\vector(0,-1){20}}
\put(300,50){\vector(0,-1){16}}
\put(170,40){\makebox(0,0){{\footnotesize $\partial _{1,n}^{S}\otimes
1_{F_{0}}$}}} 
\put(330,40){\makebox(0,0){{\footnotesize $\partial _{1,n}^{F}$}}}

\put(200,90){\vector(0,-1){20}}
\put(300,90){\vector(0,-1){20}}
\put(170,80){\makebox(0,0){{\footnotesize $\partial _{2,n}^{S}\otimes
1_{F_{0}}$}}} 
\put(330,80){\makebox(0,0){{\footnotesize $\partial _{2,n}^{F}$}}}

\put(250,110){\makebox(0,0){\mbox{\footnotesize $1\otimes \gamma
_{n-2}$}}} 
\put(250,70){\makebox(0,0){\mbox{\footnotesize $1\otimes
\gamma _{n-1}$}}} 
\put(250,10){\makebox(0,0){\mbox{\footnotesize $\gamma _{n}$}}}
\end{picture}

\noindent By Theorem 2.4 in \cite{planas2}, the middle column is exact.
Thus ${\rm ker}(\partial _{1,n}^{S}\otimes 1_{F_{0}})={\rm
im}(\partial _{2,n}^{S}\otimes 1_{F_{0}})$. Hence, $(1\otimes \gamma
_{n-1})({\rm ker}(\partial _{1,n}^{S}\otimes 1_{F_{0}}))={\rm
im}((1\otimes \gamma _{n-1})\circ (\partial _{2,n}^{S}\otimes 1_{F_{0}})
)={\rm im}(\partial _{2,n}^{F})$. Using the snake lemma, we
conclude that $E(F)_{n}={\rm ker}(\partial _{1,n}^{F})/{\rm
im}(\partial _{2,n}^{F})$.   \fidemo

\medskip

\begin{observacio}\label{prouper}
{\rm As a corollary of Proposition \ref{febrer}
we have (see also 3.1, 3.2 and 3.3 in \cite{planas2}):
\begin{itemize}
\item[$(1)$] Let $U$ be a cyclic standard $A$-algebra generated by a
degree one form $x\in U_{1}$. If $F$ is a standard $U$-module, then
$E(F)_{n}=(0:x)\cap F_{n-1}$ and ${\rm rt}(F)={\rm min}\{ r\geq 1\mid
(0:_{F}x^{r+1})=(0:_{F}x^{r})\}$. If $U=\rees{I}$ is the Rees algebra
of a principal ideal $I=(x)$ of $A$ and $F=\reesw{I}{M}$ is the Rees
module of $I$ with respect to a module $M$, then $E(I;M)=(0:x)\cap
I^{n-1}M$ and ${\rm rt}(I;M)={\rm min}\{ r\geq 1\mid
(0:_{M}x^{r+1})=(0:_{M}x^{r})\}$.
\item[$(2)$] If $\varphi :A\rightarrow B$ is a homomorphism of rings,
$U$ is standard $A$-algebra and $F$ is a standard $U$-module, then
$U\otimes B$ is a standard $B$-algebra and $F\otimes B$ is a standard
$U\otimes B$-module. Moreover ${\rm rt}_{U\otimes B\mbox{\footnotesize
-mod}}(F\otimes B)\leq {\rm rt}_{U\mbox{\footnotesize -mod}}(F)$. If
$\varphi $ is flat, $E_{U\otimes B\mbox{\footnotesize -mod}}(F\otimes
B)=E_{U\mbox{\footnotesize -mod}}(F)\otimes B$. In particular, ${\rm
rt}(F)={\rm sup}\{ {\rm rt}(F_{\mathfrak{p}})\mid \mathfrak{p}\in {\rm
Spec}(A)\} ={\rm sup}\{ {\rm rt}(F_{\mathfrak{m}})\mid \mathfrak{m}\in
{\rm Max}(A)\} $.
\item[$(3)$] If $U$ is a standard $A$-algebra, $F$ is a standard
$U$-module and $J\subseteq {\rm Ann}_{A}(F_{0})$, then $U\otimes A/J$
is a standard $A/J$-algebra, $F\otimes A/J=F$ is a standard $U\otimes
A/J$-module, $E_{U\mbox{\footnotesize -mod}}(F)_{n}=E_{U\otimes
A/J\mbox{\footnotesize -mod}}(F)_{n}$ and ${\rm
rt}_{U\mbox{\footnotesize -mod}}(F)={\rm rt}_{U\otimes
A/J\mbox{\footnotesize -mod}}(F)$.
\item[$(4)$] If ${\rm rt}(I;M)<\infty $ (for instance, if $A$
is noetherian and $M$ is a finitely generated $A$-module), then ${\rm
rt}(I;M)={\rm rt}(\agrw{I}{M})$. In particular, if $J\subset I$, then
${\rm rt}(\reesw{I}{M}\otimes A/J)={\rm rt}(I;M)$.
\end{itemize}
}\end{observacio}

\section{Proof of Theorem 2}\label{rtim}

\begin{definicio}{\rm 
Let ${\rm grt}(M)={\rm sup}\{ {\rm rt}(I;M)\mid I \mbox{ ideal of }
A\} $ denote the supremum (possibly infinite) of all relation types of
ideals $I$ of $A$ with respect to the $A$-module $M$, and let us call
it the {\em global relation type} of $M$. Remark that:
\begin{itemize}
\item[$(1)$] If $M=A$, ${\rm grt}(A)={\rm sup}\{ {\rm rt}(I)\mid I
\mbox{ ideal of } A\} $. We will prove that for an excellent ring $A$,
${\rm grt}(A)<\infty $ is equivalent to ${\rm dim}\, A\leq 1$.
\item[$(2)$] Since ${\rm rt}(I;M)={\rm sup}\{ {\rm
rt}(I_{\mathfrak{p}};M_{\mathfrak{p}})\mid \mathfrak{p}\in {\rm
Spec}(A)\} = {\rm sup}\{ {\rm
rt}(I_{\mathfrak{m}};M_{\mathfrak{m}})\mid \mathfrak{m}\in {\rm
Max}(A)\}$, then ${\rm grt}(M)={\rm sup}\{ {\rm
grt}(M_{\mathfrak{p}})\mid \mathfrak{p}\in {\rm Spec}(A)\} = {\rm
sup}\{ {\rm grt}(M_{\mathfrak{m}})\mid \mathfrak{m}\in {\rm
Max}(A)\}$.
\item[$(3)$] If necessary to specify, we will write ${\rm grt}(M)={\rm
grt}_{A}(M)$ when considering $M$ an $A$-module. For instance, if
$J\subseteq {\rm Ann}_{A}(M)=\{ x\in A\mid xM=0\} $, then
$\reesw{(I+J)/J}{M}=\reesw{I}{M}$. Thus ${\rm rt}(I;M)={\rm
rt}((I+J)/J;M)$ and ${\rm grt}_{A}(M)={\rm grt}_{A/J}(M)$.
\end{itemize}
}\end{definicio}

\noindent {\bf Theorem 2} {\em 
Let $A$ be a commutative ring, $\mathcal{W}$ a set of ideals of $A$,
$I\in \mathcal{W}$ and $N\subseteq M$ two $A$-modules. Let
$s(N,M;\mathcal{W})$ denote the strong uniform number for the pair
$(N,M)$ with respect to $\mathcal{W}$. Then $s(N,M;\{ I\} )\leq {\rm
rt}(I;M/N)\leq {\rm max}({\rm rt}(I;M),s(N,M;\{ I\}))$. In particular,
$s(N,M;\mathcal{W})\leq {\rm sup}\{ {\rm rt}(J;M/N)\mid J\in
\mathcal{W}\}$ and $s(N,M)\leq {\rm grt}(M/N)$.}

\bigskip

\noindent {\em Proof}. Let $F=\reesw{I}{M/N}$, $G=\reesw{I}{M}$,
$H={\bf S}(I)\otimes M$, $\varphi :G\rightarrow F$ the surjective
graded morphism of standard ${\bf S}(I)$-algebras defined by $\varphi
_{n}:G_{n}=I^{n}M\rightarrow I^{n}M/I^{n}M\cap N= I^{n}M+N/N=F_{n}$
and $\gamma :H\rightarrow G$ induced by the natural graded morphism
$\alpha :{\bf S}(I)\rightarrow \rees{I}$. By Lemma \ref{gener},
$s(E(\varphi ))\leq {\rm s}(E(\varphi \circ \gamma ))\leq {\rm
max}(s(E(\varphi )),s(E(\gamma )))$. But, $s(E(\varphi \circ \gamma
))={\rm rt}(I;M/N)$ and $s(E(\gamma ))={\rm rt}(I;M)$. Finally, since
$E(\varphi )_{n}=(I^{n}M\cap N)/I(I^{n-1}M\cap N)$, then $s(E(\varphi
))= s(N,M;\{ I\} )$. In particular, $s(N,M;\{ I\} )\leq {\rm
rt}(I;M/N)\leq {\rm sup}\{ {\rm rt}(J;M/N)\mid J\in \mathcal{W}\}$, and
taking the supremum over all ideals $I$ of $\mathcal{W}$,
$s(N,M;\mathcal{W})\leq {\rm sup}\{ {\rm rt}(J;M/N)\mid J\in
\mathcal{W}\}$.  \fidemo

\begin{corollari} {\sc Artin-Rees Lemma}.
Let $A$ be a commutative ring, $I$ an ideal of $A$ and $N\subseteq M$
two $A$-modules. If ${\rm rt}(I;M/N)<\infty $ then $s(N,M;\{ I\}
)<\infty $. In particular, if $A$ is noetherian and $M$ is finitely
generated, there exists an integer $s\geq 1$ such that, for all
integers $n\geq s$, $I^{n}M\cap N=I^{n-s}(I^{s}M\cap N)$.
\end{corollari}

\begin{corollari}\label{princip}
{\sc O'Carroll \cite{ocarroll1}}. Let $A$ be a noetherian ring and let
$M$ be a finitely generated $A$-module. Then ${\rm sup}\{ {\rm
rt}((x);M)\mid x\in A\} <\infty$. In particular, if $N\subset M$,
there exists an integer $s\geq 1$ such that, for all integers $n\geq
s$ and for all $x\in A$, $x^{n}M\cap N=x^{n-s}(x^{s}M\cap N)$.
\end{corollari}

\noindent {\em Proof}. Following very closely the proof of O'Carroll
in \cite{ocarroll1}, let $0=Q_{1}\cap \ldots \cap Q_{r}$ be a minimal
primary decomposition of $0$ in $M$,
$r_{M}(Q_{i})=r(Q_{i}:M)=\mathfrak{p}_{i}\in {\rm Spec}(A)$, and let
$s\geq 1$ be an integer such that, for all $i=1,\ldots ,r$,
$\mathfrak{p}_{i}^{s}M\subseteq Q_{i}$. Then, for all $x\in A$, ${\rm
rt}((x);M)\leq s$. Indeed, if $x\in \mathfrak{p}_{i}$, $x^{n+s}\in
\mathfrak{p}_{i}^{n+s}$ and $(Q_{i}:x^{n+s})=M$. If $x\not\in
\mathfrak{p}_{i}$, $x^{n+s}\not\in \mathfrak{p}_{i}^{n+s}$ and
$(Q_{i}:x^{n+s})=Q_{i}$. Therefore, for all $n\geq 0$,
$(0:x^{n+s})=(\cap _{i}Q_{i}:x^{n+s})= \cap _{i}(Q_{i}:x^{n+s})=\cap
_{x\not\in \mathfrak{p}_{i}}Q_{i}$. In particular,
$(0:x^{s+1})=(0:x^{s})$ and ${\rm rt}((x);M)\leq s$. We finish by
applying Theorem 2. \fidemo

\begin{observacio}{\rm 
Let $A$ be a noetherian ring and let $M$ be a finitely generated
$A$-module. Let ${\rm grt}^{i}(M)={\rm sup}\{ {\rm rt}(I;M)\mid
\mu (I)\leq i\}$. By \ref{princip},
${\rm grt}^{1}(M)<\infty $. Using the example of Wang \cite{wang1} and
Theorem 2, we know ${\rm grt}^{3}(M)$ might be
infinite. We do not know whether ${\rm grt}^{2}(M)$ is finite.
}\end{observacio}

\section{Rings of finite global relation type have dimension
one}\label{1dim} 

\begin{observacio}\label{rest05} {\rm Let $A$ be a commutative
ring and let $r\geq 1$ denote an integer. Consider the following
conditions:
\begin{itemize}
\item[$(a)$] ${\rm rt}(I)\leq r$ for every three-generated ideal $I$
of $A$.
\item[$(b)$] $E(I)_{r+1}=0$ for every three-generated ideal $I$ of $A$.
\item[$(c)$] $(x,y)(x,y,z)^{r}:z^{r+1}=(x,y)(x,y,z)^{r-1}:z^{r}$ for
all $x,y,z\in A$.
\item[$(d)$] $(x^{r}y)^{r}\in
(x^{r+1},y^{r+1})(x^{r+1},y^{r+1},x^{r}y)^{r-1}$ for all $x,y\in A$.
\end{itemize}
Then $(a)\Rightarrow (b)\Rightarrow (c)\Rightarrow (d)$.
}\end{observacio}

\noindent {\em Proof}. Implication $(a)\Rightarrow (b)$ follows from
the definitions. Implication $(b)\Rightarrow (c)$ holds in general: if
$I$ is generated by $x_{1},\ldots ,x_{d}$ and if $E(I)_{n}=0$, then
$(x_{1},\ldots ,x_{d-1})I^{n-1}:x_{d}^{n} =(x_{1},\ldots
,x_{d-1})I^{n-2}:x_{d}^{n-1}$ (Lemma 4.2 \cite{planas2}). Finally,
$(d)$ follows from $(c)$ taken $x,y,z\in A$ as
$x^{r+1},y^{r+1},x^{r}y$.  \fidemo

\medskip

In order to prove $(d)\Rightarrow {\rm dim}\, A\leq 1$ let us recall
some definitions. A set of elements $x_{1},\ldots ,x_{m}$ of an ideal
$J$ of $A$ are called $J$-{\em independent} if every form in
$A[T_{1},\ldots ,T_{m}]$ vanishing at $x_{1},\ldots ,x_{m}$ has all
its coefficients in $J$. If $I=(x_{1},\ldots ,x_{m})$ and $I\subset
J$, then $x_{1},\ldots ,x_{m}$ are $J$-independent if and only if the
natural graded morphism of standard $(A/J)$-algebras
$(A/J)[X_{1},\ldots ,X_{m}]\rightarrow \rees{I}\otimes (A/J)$ is an
isomorphism ($X_{1},\ldots ,X_{m}$ algebraically independent over
$A/J$). If $(A,\mathfrak{m})$ is noetherian local, then the maximum
number of $\mathfrak{m}$-independent elements in $\mathfrak{m}$ is
equal to ${\rm dim}\, A$ \cite{valla}.

\begin{proposicio}\label{rest1} Let $A$ be a
noetherian ring. If there exists an integer $r\geq 1$ such that
$(x^{r}y)^{r}\in (x^{r+1},y^{r+1})(x^{r+1},y^{r+1},x^{r}y)^{r-1}$ for
all $x,y\in A$, then ${\rm dim}\, A\leq 1$.
\end{proposicio}

\noindent {\em Proof}. Since the hypothesis localizes, we may assume
$(A,\mathfrak{m},k)$ is a noetherian local ring. Suppose ${\rm
dim}\, A\geq 2$. Then there exists two
$\mathfrak{m}$-independent elements $x,y$. In particular, if
$I=(x,y)$, $\overline{\alpha}:k[X,Y]\rightarrow
\rees{I}/\mathfrak{m}\rees{I}$ defined by
$\overline{\alpha}(X)=x+\mathfrak{m}I^{2}$,
$\overline{\alpha}(Y)=y+\mathfrak{m}I^{2}$ is a graded isomorphism of
standard $k$-algebras. By hypothesis $(x^{r}y)^{r}\in
(x^{r+1},y^{r+1})(x^{r+1},y^{r+1},x^{r}y)^{r-1}$ which is generated by
the elements $x^{(i+1)(r+1)+lr}y^{j(r+1)+l},
x^{i(r+1)+lr}y^{(j+1)(r+1)+l}$, $i,j,l\geq 0$, $i+j+l=r-1$.  The
$k$-vector space isomorphism $\overline{\alpha}_{r(r+1)}$ assures the
membership of $(X^{r}Y)^{r}$ into the $k$-vector space spanned by the
elements $X^{(i+1)(r+1)+lr}Y^{j(r+1)+l}$,
$X^{i(r+1)+lr}Y^{(j+1)(r+1)+l}$, $i,j,l\geq 0$, $i+j+l=r-1$. Since all
of them are elements of a $k$-basis of $k[X,Y]_{r(r+1)}$, then either
$(X^{r}Y)^{r}=X^{(i+1)(r+1)+lr}Y^{j(r+1)+l}$ or either $(X^{r}Y)^{r}=
X^{i(r+1)+lr}Y^{(j+1)(r+1)+l}$, for some $i,j,l\geq 0$, $i+j+l=r-1$.
But, it is not difficult to see that there are not integers $i,j,l\geq
0$ verifying any of both equations. \fidemo

\begin{observacio}{\rm 
The underlying idea in the proof of Proposition \ref{rest1} is that
for any two $\mathfrak{m}$-independent elements $x,y$ of $A$, do not
exist $r$-relations
$T_{1}f(T_{1},T_{2},T_{3})+T_{2}g(T_{1},T_{2},T_{3})-T_{3}^{r}$,
$f,g$ forms of degree $r-1$, among the three ordered elements
$x^{r+1},y^{r+1},x^{r}y$. In
particular, $T_{1}^{r}T_{2}-T_{3}^{r+1}$ must be an effective
$(r+1)$-relation among the three ordered elements $x^{r+1},y^{r+1},x^{r}y$
(since any form of degree $r$ dividing
$T_{1}^{r}T_{2}-T_{3}^{r+1}$ should contain $T_{3}^{r}$
as an additive factor). }\end{observacio}

\begin{observacio}{\rm 
There exist (necessarily non noetherian) local rings with ${\rm
grt}(A)<\infty $, but ${\rm dim}\, A\geq 2$. For example, a valuation
ring $A$ is Pr\"ufer, thus ${\rm grt}(A)=1$ \cite{costa},
\cite{planas1}, but its dimension is not necessarily 1 or less.
}\end{observacio}

\section{Artinian modules have finite global relation type}\label{0dim}

\begin{proposicio}\label{artim} Let $A$ be a commutative ring,
$I$ an ideal of $A$ and $M$ an $A$-module. If $I^{s}M=0$ for some
$s\geq 1$, then ${\rm rt}(I;M)\leq s$. If ${\rm rt}(I)=1$ and $I$ is
finitely generated, then $I\neq 0$ if and only if $I^{s}\neq 0$ for
all $s\geq 1$.  If $(A,\mathfrak{m})$ is artinian local, then ${\rm
grt}(M)<\infty $ and ${\rm grt}(A)=1$ if and only if $A$ is a
field. If $A$ is an artinian ring, then ${\rm grt}(M)<\infty $.
\end{proposicio}

\noindent {\em Proof}. If $I^{s}M=0$ and $\varphi :G\rightarrow
\reesw{I}{M}$ is a symmetric presentation of $\reesw{I}{M}$, then
${\rm ker}\varphi _{n}=G_{n}$ for all $n\geq s$. Thus, for all $n\geq
s+1$, $E(I)_{n}=G_{n}/V_{1}G_{n-1}=0$ and ${\rm rt}(I;M)\leq s$.  In
order to prove the second assertion we may suppose that
$(A,\mathfrak{m},k)$ is local. If ${\rm rt}(I)=1$, the natural graded
morphism of standard $k$-algebras ${\bf
S}^{k}(I/\mathfrak{m}I)\rightarrow \rees{I}/\mathfrak{m}\rees{I}$ is
an isomorphism. If $I\neq 0$, then ${\bf S}^{k}(I/\mathfrak{m}I)$ is a
polynomial ring in $\mu (I)$ variables, thus
$I^{s}/\mathfrak{m}I^{s}\neq 0$ for all $s\geq 1$ and $I^{s}\neq 0$
since $I$ is finitely generated. If $(A,\mathfrak{m})$ is artinian
local, there exists an integer $s\geq 1$, such that $I^{s}M\subseteq
\mathfrak{m}^{s}M=0$ for every ideal $I$ of $A$. Thus ${\rm
grt}(M)\leq s$. Moreover, if ${\rm grt}(A)=1$, then ${\rm
rt}(\mathfrak{m})=1$ and $\mathfrak{m}^{s}=0$. Hence $\mathfrak{m}=0$
and $A$ is a field. If $A$ is artinian, it has a finite number of
maximal ideals and since ${\rm grt}(M)={\rm sup}\{ {\rm
grt}(M_{\mathfrak{m}}) \mid \mathfrak{m}\in {\rm Max}(A)\} $, then
${\rm grt}(M)<\infty $. \fidemo

\begin{observacio}{\rm The minimum integer $s\geq 1$ such that $I^{s}=0$,
for a nilpotent ideal $I$, is not necessarily equal to its relation
type. For example, take $I=(x,y)\subset A= k\lbrack\!\lbrack
X,Y\rbrack\!\rbrack /(X^{n},Y^{n})$, where $x,y$ denote the classes of
$X,Y$ in $A$. Then $I^{2n-1}=0$, $I^{2n-2}\neq 0$ and ${\rm
rt}(I)=n$. Indeed, since $yI^{n-2}:x^{n-1}\varsubsetneq
yI^{n-1}:x^{n}=A$, then $E(I)_{n}\neq 0$ and ${\rm rt}(I)\geq
n$. Moreover, since $(0:y)\cap I^{p-1}=(x^{p-n}y^{n-1})= x((0:y)\cap
I^{p-2})$ for all $p\geq n+1$, then $E(I)_{p}=0$ for all $p\geq n+1$
and ${\rm rt}(I)\leq n$ (Proposition 4.5 \cite{planas2}). 
}\end{observacio}

\begin{observacio}{\rm There exist (necessarily non noetherian) local
rings with ${\rm dim}\, A=0$, but ${\rm grt}(A)=\infty $. For example,
$A=k[T_{1},\ldots ,T_{m},\ldots ]/(T_{1}^{2},\ldots
,T_{m}^{m+1},\ldots )$, with $k$ a field, is a zero dimensional local
ring.  If $t_{m}$ denotes the residue class of $T_{m}$,
$(0:t_{m}^{m})\varsubsetneq (0:t_{m}^{m+1})=A$ and ${\rm
rt}((t_{m}))=m+1$. }\end{observacio}

\section{Proof of Theorem 3 in the local case}\label{jor}

We first need to reduce to Cohen-Macaulay local rings  and modules.

\begin{lema}\label{pasquom}
Let $A$ be a noetherian ring, $I$ an ideal of $A$ and $N\subseteq M$
two finitely generated $A$-modules such that $I^{t}N=0$ for a certain
integer $t\geq 1$. Then ${\rm rt}(I;M)\leq {\rm rt}(I;M/N)+t$. In
particular, if $I$, $J$ are two ideals of $A$ such that $I^{t}J=0$ for
a certain integer $t\geq 1$, then ${\rm rt}(I)\leq {\rm
rt}((I+J)/J)+t$.
\end{lema}

\noindent {\em Proof}. Let $s={\rm s}(N,M;\{ I\})$ be the strong
uniform number for the pair $(N,M)$ with respect to the set of ideals
$\{ I\} $. If $n\geq s+t$, then $I^{n}M\cap N=I^{n-s}(I^{s}M\cap
N)\subseteq I^{n-s}N\subseteq I^{t}N=0$. Let $F=\reesw{I}{M/N}$,
$G=\reesw{I}{M}$ and $\varphi:G\rightarrow F$ defined by $\varphi
_{n}:G_{n}=I^{n}M\rightarrow I^{n}M/I^{n}M\cap N=I^{n}M+N/N=F_{n}$.
We have ${\rm ker}\varphi _{n}=I^{n}M\cap N=0$ for all $n\geq
s+t$. Therefore, using Remark \ref{prougranm} and Theorem 2, ${\rm
rt}(I;M)={\rm rt}(G)\leq {\rm max}({\rm rt}(F),s+t)\leq {\rm max}({\rm
rt}(I;M/N),{\rm rt}(I;M/N)+t)={\rm rt}(I;M/N)+t$.  \fidemo

\begin{corollari}\label{longfinm}
Let $(A,\mathfrak{m})$ be a noetherian local ring, $M$ a finitely
generated $A$-module and $N\subseteq M$ a submodule of finite
length. Then ${\rm grt}(M)\leq {\rm grt}(M/N)+{\rm length}(N)$.
\end{corollari}

\noindent {\em Proof}. If ${\rm length}(N)=t$, then $I^{t}N\subset
\mathfrak{m}^{t}N=0$ for every ideal $I$ of $A$. Thus, by Lemma
\ref{pasquom}, ${\rm rt}(I;M)\leq {\rm rt}(I;M/N)+t\leq {\rm
grt}(M/N)+t$. Taking the supremum over all ideals $I$ of $A$, we have
${\rm grt}(M)\leq {\rm grt}(M/N)+t$.  \fidemo

\medskip

\noindent Next lemma is a generalization to modules of a well known
result for rings (see, for instance, \cite{trung2}).

\begin{lema}\label{mprim}
Let $(A,\mathfrak{m})$ be a one dimensional Cohen-Macaulay local
ring. Let $M$ be a maximal Cohen-Macaulay
module. If $I$ is an $\mathfrak{m}$-primary ideal of $A$, ${\rm
rt}(I;M)\leq e(A)$.
\end{lema}

\noindent {\em Proof}. Applying the tensor product $-\otimes
A[t]_{\mathfrak{m}[t]}$, we may assume that the residue field
$k=A/\mathfrak{m}$ is infinite. By Theorem 1.1 in \cite{sally}, $\mu
(I)\leq e(A)=e$ and $\mu (I^{e})\leq e< {e+1\choose 1} $. By Theorem
2.3 in \cite{sally}, there exists $y_{0}\in I$ such that
$I^{e}=y_{0}I^{e-1}$. In particular, for all $n\geq e$,
$I^{n}=y_{0}I^{n-1}$. Since $\mathfrak{m}\not\subset Z(M)$, then
$y_{0}\not\in Z(M)$. Consider the complex of $A$-modules:
\begin{eqnarray*}
{\bf \Lambda}_{2}(I)\otimes I^{n-2}M\buildrel
\partial _{2,n}\over \longrightarrow I\otimes I^{n-1}M\buildrel
\partial _{1,n}\over \longrightarrow I^{n}M\longrightarrow 0\, ,
\end{eqnarray*}
where $\partial _{2,n}((x\wedge y)\otimes z)=y\otimes xz-x\otimes yz$
and $\partial _{1,n}(x\otimes t)=xt$, $x,y\in I$, $z\in I^{n-2}M$ and
$t\in I^{n-1}M$. By Proposition \ref{febrer}, $E(I;M)_{n}={\rm
ker}\partial _{1,n}/{\rm im}\partial _{2,n}$. Let us see ${\rm
ker}\partial _{1,n}={\rm im}\partial _{2,n}$ for all $n\geq
e+1$. Indeed, take $u=\sum x_{i}\otimes y_{0}z_{i}\in {\rm
ker}\partial _{1,n}$, $x_{i}\in I$, $z_{i}\in I^{n-2}M$. Then
$0=\partial _{1,n}(u)=y_{0}\sum x_{i}z_{i}$ and, since $y_{0}\not\in
Z(M)$, $\sum x_{i}z_{i}=0$. Therefore, if $v=\sum (y_{0}\wedge
x_{i})\otimes z_{i}\in {\bf \Lambda}_{2}(I)\otimes I^{n-2}M$, then
$\partial _{2,n}(v)=\sum x_{i}\otimes y_{0}z_{i}-\sum y_{0}\otimes
x_{i}z_{i}=u$. So $E(I;M)_{n}=0$ for all $n\geq e+1$ and ${\rm
rt}(I;M)\leq e(A)$. \fidemo

\begin{notacions}\label{conv} {\rm Let $(A,\mathfrak{m})$ be a one 
dimensional noetherian local ring. Denote by $\mathfrak{q}_{1},\ldots
,\mathfrak{q}_{s}$ the minimal primary components of $(0)$. If $A$ is
Cohen-Macaulay, $(0)=\mathfrak{q}_{1}\cap \ldots \cap
\mathfrak{q}_{s}$ is a minimal primary decomposition of $(0)$. If $A$
is not Cohen-Macaulay, there exist an $\mathfrak{m}$-primary ideal
$\mathfrak{q}_{s+1}$ such that $(0)=\mathfrak{q}_{1}\cap \ldots \cap
\mathfrak{q}_{s}\cap \mathfrak{q}_{s+1}$ is a minimal primary
decomposition of $(0)$. Let $n\geq 1$ be the minimum integer such that
$\mathfrak{n}(A)^{n}=0$. Let $n_{i}\geq 1$ be the minimum integer such
that $\mathfrak{p}_{i}^{n_{i}}\subset \mathfrak{q}_{i}$,
$\mathfrak{p}_{i}=r(\mathfrak{q}_{i})$, $i=1,\ldots ,s$.  
For each $1\leq i_{1}<\ldots <i_{l}\leq s$, denote
$t_{i_{1},\ldots ,i_{l}}={\rm max}\{ n_{i}\mid i\neq i_{1},\ldots
,i_{l}\} $ and $e(A)$ the multiplicity of $A$. Finally, set ${\rm
brt}(A)={\rm max}\{ n, e(A/(\mathfrak{q}_{i_{1}}\cap \ldots \cap
\mathfrak{q}_{i_{l}}))+ t_{i_{1},\ldots ,i_{l}} \mid 1\leq
i_{1}<\ldots <i_{l}\leq s \}$, which is finite. 
}\end{notacions}

\begin{proposicio}\label{bo} 
Let $(A,\mathfrak{m})$ be a one dimensional noetherian local ring and
$J={\rm H}^{0}_{\mathfrak{m}}(A)$. Let $M$ be a one dimensional
finitely generated $A$-module and $N={\rm H}^{0}_{\mathfrak{m}}(M)$.
Then ${\rm grt}(A)\leq {\rm brt}(A/J)+{\rm length}(J)$ and ${\rm
grt}(M)\leq {\rm brt}(A/J)+{\rm length}(N)$. If $A$ and $M$ are
Cohen-Macaulay, ${\rm grt}(M)\leq {\rm brt}(A)$.
\end{proposicio}

\noindent {\em Proof}. Since ${\rm length}(N)=t<\infty$, by Corollary
\ref{longfinm}, ${\rm grt}(M)\leq {\rm grt}(M/N)+t$. Since
$JM\subseteq N$, then $J\subseteq {\rm Ann}_{A}(M/N)$ and ${\rm
grt}_{A}(M/N)={\rm grt}_{A/J}(M/N)$. We thus may assume $A$ is a one
dimensional Cohen-Macaulay ring and $M$ is a maximal Cohen-Macaulay
module. Let us prove ${\rm grt}(M)\leq {\rm brt}(A)$.  If $I\subset
\mathfrak{n}(A)$, $I^{n}\subset \mathfrak{n}(A)^{n}=0$ and ${\rm
rt}(I;M)\leq n$ (Proposition \ref{artim}). If $I\not\subset
\mathfrak{N}(A)$, let $1\leq i_{1}<\ldots <i_{l}\leq s$ be all the
subindexes $i_{j}$ such that $I\not\subset \mathfrak{p}_{i_{j}}$. Set
$J_{i_{1},\ldots ,i_{l}}=\mathfrak{q}_{i_{1}}\cap \ldots \cap
\mathfrak{q}_{i_{l}}$. Then $I^{t_{i_{1},\ldots
,i_{l}}}J_{i_{1},\ldots ,i_{l}}\subseteq \mathfrak{q}_{1}\cap \ldots
\cap \mathfrak{q}_{s}=0$ and, by Lemma \ref{pasquom}, ${\rm
rt}(I;M)\leq {\rm rt}(I;M/J_{i_{1},\ldots ,i_{l}}M)+t_{i_{1},\ldots
,i_{l}}={\rm rt}((I+J_{i_{1},\ldots ,i_{l}})/J_{i_{1},\ldots
,i_{l}};M/J_{i_{1},\ldots ,i_{l}}M)+t_{i_{1},\ldots ,i_{l}}$. But,
$(I+J_{i_{1},\ldots ,i_{l}})/J_{i_{1},\ldots ,i_{l}}$ is an
$\mathfrak{m}/J_{i_{1},\ldots ,i_{l}}$-primary ideal of the one
dimensional Cohen-Macaulay local ring $A/J_{i_{1},\ldots ,i_{l}}$ and
$M/J_{i_{1},\ldots ,i_{l}}M$ is a maximal Cohen-Macaulay module.
Therefore, by Lemma \ref{mprim}, ${\rm rt}((I+J_{i_{1},\ldots
,i_{l}})/J_{i_{1},\ldots ,i_{l}};M/J_{i_{1},\ldots ,i_{l}}M)\leq
e(A/J_{i_{1},\ldots ,i_{l}})$. \fidemo

\begin{exemple}\label{doman} {\rm
Let $(A,\mathfrak{m})$ be a one dimensional noetherian local ring. If
$A$ is reduced, then ${\rm grt}(A)\leq e(A)+1$. If $A$ is a domain,
then ${\rm grt}(A)\leq e(A)$.  }\end{exemple}

\noindent {\em Proof}. Since $A$ is Cohen-Macaulay, by Proposition
\ref{bo}, ${\rm grt}(A)\leq {\rm brt}(A)$.  Following the notations in
\ref{conv}, if $A$ is reduced, $n=n_{1}=\ldots =n_{s}=1$,
$t_{i_{1},\ldots ,i_{l}}=1$ for all $(i_{1},\ldots ,i_{l})\neq
(1,\ldots ,s)$ and $t_{1,\ldots ,s}=0$. Since $e(A/J)\leq e(A)$, then
${\rm brt}(A)\leq e(A)+1$. If $A$ is a domain, then $n=1$, $n_{1}=1$,
$t_{1}=0$ and ${\rm brt}(A)=e(A)$. \fidemo

\begin{exemple}\label{localtg} {\rm
Let $k$ be a field and $g\geq 1$ an integer. Set
$R=k[t^{g+1},t^{g+2},\ldots ,t^{2g+1}]\subset k[t]$ ($t$ a variable
over $k$), $\mathfrak{n}=(t^{g+1},t^{g+2},\ldots ,t^{2g+1})$,
$A=R_{\mathfrak{n}}$ and $\mathfrak{m}=\mathfrak{n}R_{\mathfrak{n}}$.
Then $(A,\mathfrak{m},k)$ is a one dimensional notherian local domain
and ${\rm grt}(A)=e(A)=g+1$.  }\end{exemple}

\noindent {\em Proof}. By Example \ref{doman}, ${\rm grt}(A)\leq
e(A)$.  For all $n\geq 1$, $\mathfrak{m}^{n}=(t^{(g+1)n},\ldots
,t^{(g+1)n+g})$, $\mu (\mathfrak{m}^{n})=g+1$ and $e(A)=g+1$. For
$n\geq 2$, take $I=(t^{g+1},t^{g+2})$ and
$J=_{g,n-1}=t^{g+1}I^{n-2}:t^{(g+2)(n-1)}$. Remark that
$J_{g,n-1}\subseteq J_{g,n}$ and that $E(I)_{n}=0$ if and only if
$J_{g,n-1}=J_{g,n}$ (Proposition 4.5 \cite{planas2}). If $g=1$, then
$I=\mathfrak{m}$, $E(I)_{2}\neq 0$ and $2\leq {\rm rt}(I)\leq
e(A)=2$. Suppose $g\geq 2$. Then
$J_{g,1}=t^{g+1}:t^{g+2}=\mathfrak{m}$.  Moreover, $t^{(g+2)g}\not\in
t^{g+1}I^{g-1}$ and $\mathfrak{m}\subseteq J_{g,g}\varsubsetneq
A$. Moreover $J_{g,g+1}=A$. Thus, $E(I)_{n}=0$ for all $2\leq n\leq g$
and $E(I)_{g+1}\neq 0$. Hence $g+1\leq {\rm rt}(I)\leq e(A)=g+1$,
${\rm rt}(I)=g+1$ and ${\rm grt}(A)=g+1$. Remark that
$\mathfrak{m}^{n}=t^{g+1}\mathfrak{m}^{n-1}$ for all $n\geq 2$. So the
reduction number of $\mathfrak{m}$ is ${\rm rn}(\mathfrak{m})=1$ and
$1<{\rm rt}(\mathfrak{m})\leq {\rm rn}(\mathfrak{m})+1=2$ \cite{trung2}
while ${\rm grt}(A)=g+1$.  \fidemo

\begin{exemple}\label{eigrt}{\rm
Let $k$ be a field, $a\geq 1$ a positive integer and
$A=k\lbrack\!\lbrack X, Y\rbrack\!\rbrack /(X^{a}Y)$. Then $A$ is a
one dimensional complete intersection local ring with ${\rm
grt}(A)={\rm brt}(A)=a+1$.  }\end{exemple}

\noindent {\em Proof}. By Proposition \ref{bo}, ${\rm grt}(A)\leq {\rm
brt}(A)$. Let $x,y$ denote the residue classes of $X,Y$ and let
$\mathfrak{m}=(x,y)$ be the maximal ideal of $A$. Since $\mu
(\mathfrak{m}^{n})=a+1$ for all $n\geq a$, the multiplicity of $A$ is
$e(A)=a+1$. The minimal primary decomposition of $A$ is
$(0)=\mathfrak{q}_{1}\cap \mathfrak{q}_{2}$,
$\mathfrak{q}_{1}=(x^{a})$, $\mathfrak{q}_{2}=(y)$. Following the
notations in \ref{conv}, $\mathfrak{p}_{1}=(x)$,
$\mathfrak{p}_{2}=(y)$, $\mathfrak{n}(A)=(xy)$, $n=n_{1}=a$,
$n_{2}=1$, $t_{1}=n_{2}=1$, $t_{2}=n_{1}=a$, $t_{1,2}=0$. Moreover,
$A/\mathfrak{q}_{1}= k\lbrack\!\lbrack X, Y\rbrack\!\rbrack /(X^{a})$
and $e(A/\mathfrak{q}_{1})=a$; $A/\mathfrak{q}_{2}= k\lbrack\!\lbrack
X\rbrack\!\rbrack $ and $e(A/\mathfrak{q}_{2})=1$. Therefore, ${\rm
brt}(A)=a+1$. On the other hand, $x((0:y)\cap
\mathfrak{m}^{a-1})=(x^{a+1})\varsubsetneq (x^{a})=(0:y)\cap
\mathfrak{m}^{a}$. Thus $E(\mathfrak{m})_{a+1}\neq 0$ and ${\rm
rt}(\mathfrak{m})\geq a+1$ (Proposition 4.5 \cite{planas2}). \fidemo

\section{Final proofs}\label{final}

\begin{lema}\label{quiqm} {\rm Let $(A,\mathfrak{m})$ be a 
one dimensional Cohen-Macaulay local ring with a unique minimal prime
$\mathfrak{p}$ and let $n\geq 1$ be such that $\mathfrak{p}^{n}=0$. If
$M$ is a maximal Cohen-Macaulay $A$-module, then ${\rm grt}(M)\leq
{\rm max}\{ n,e(A)\} ={\rm brt}(A)$. Moreover, if $A/\mathfrak{p}$ is
a discrete valuation ring, then ${\rm grt}(M)\leq {\rm max}\{ n,\sum
_{i=0}^{n-1}\mu (\mathfrak{p}^{i})\} $.  }\end{lema}

\noindent {\em Proof}. By Proposition \ref{bo}, ${\rm grt}(M)\leq {\rm
brt}(A)$.  If $I\subseteq \mathfrak{p}$, then $I^{n}\subseteq
\mathfrak{p}^{n}=0$ and ${\rm rt}(I;M)\leq n$. If $I\not\subset
\mathfrak{p}$, then $I$ is an $\mathfrak{m}$-primary ideal of a one
dimensional Cohen-Macaulay local ring. Hence, by Lemma \ref{mprim},
${\rm rt}(I;M)\leq e(A)$. Remark that ${\rm brt}(A)={\rm max}\{ n,
e(A)\} $. If moreover, $A/\mathfrak{p}$ is a discrete valuation ring,
there exists $u\in A$ such that $\mathfrak{m}=uA+\mathfrak{p}$. Thus,
for $r\geq n$, $\mathfrak{m}^{r}=\sum
_{i=0}^{n-1}u^{r-i}\mathfrak{p}^{i}$ and for $r\gg 1$, $e(A)=\mu
(\mathfrak{m}^{r})=\mu (\sum _{i=0}^{n-1}u^{r-i}\mathfrak{p}^{i})\leq
\sum _{i=0}^{n-1}\mu (\mathfrak{p}^{i})$. \fidemo

\begin{exemple}{\rm Let $k$ be a field, $a,b\geq 1$ 
two positive integers and $A=k\lbrack\!\lbrack X, Y\rbrack\!\rbrack
/(X^{a},X^{b}Y)$. Then $A$ is a one dimensional noetherian local ring
with ${\rm grt}(A)=a$. Moreover, if $a\leq b$, then $J={\rm
H}^{0}_{\mathfrak{m}}(A)=0$ and ${\rm brt}(A)=a$. If $a>b$, then
$J={\rm H}^{0}_{\mathfrak{m}}(A)\neq 0$ and ${\rm brt}(A/J)+{\rm
length}(J)=a+b$. }\end{exemple}

\noindent {\em Proof}. Let $x,y$ denote the residue classes of $X,Y$
and let $\mathfrak{m}=(x,y)$ be the maximal ideal of $A$. Remark that
${\rm rt}((x))=a\leq {\rm grt}(A)$.  If $a\leq b$,
$A=k\lbrack\!\lbrack X, Y\rbrack\!\rbrack /(X^{a})$ is a one
dimensional Cohen-Macaulay ring with the unique minimal prime
$(x)$. By Lemma \ref{quiqm}, $a\leq {\rm grt}(A)\leq {\rm
brt}(A)={\rm max}\{ a, e(A)\} =a$. If $a>b$, then $I(x^{a-1})\subseteq
(x,y)(x^{a-1})\subseteq (x^{a},x^{b}y)=0$. By Lemma
\ref{pasquom}, ${\rm rt}(I)\leq {\rm
rt}((I+(x^{a-1}))/(x^{a-1}))+1\leq {\rm grt}(A/(x^{a-1}))+1$. But,
$A/(x^{a-1})=k\lbrack\!\lbrack X, Y\rbrack\!\rbrack
/(X^{a-1},X^{b}Y)$. Repeating the same argument, we get $a\leq {\rm
grt}(A)\leq {\rm grt}(A/(x^{a-(a-b)}))+(a-b)={\rm
grt}(k\lbrack\!\lbrack X, Y\rbrack\!\rbrack
/(X^{b}))+(a-b)=b+(a-b)=a$. On the other hand, $J={\rm
H}^{0}_{\mathfrak{m}}(A)=(0:\mathfrak{m}^{a-b})=(x^{b})$ and ${\rm
lenght}(J)=a-b$. $A/J=k\lbrack\!\lbrack X, Y\rbrack\!\rbrack /(X^{b})$
and, as before, ${\rm brt}(A/J)=b$. Thus, ${\rm brt}(A/J)+{\rm
length}(J)=a+b$. \fidemo

\bigskip

\noindent {\bf Theorem 3} {\em 
Let $A$ be an excellent (in fact $J-2$) ring. The following conditions
are equivalent:
\begin{itemize}
\item[$(i)$] ${\rm grt}(M)<\infty $ for all finitely generated
$A$-module $M$.
\item[$(ii)$] ${\rm grt}(A)<\infty $.
\item[$(iii)$] There exists an $r\geq 1$ such that ${\rm
rt}(I)\leq r$ for every three-generated ideal $I$ of $A$.
\item[$(iv)$] There exists an $r\geq 1$ such that
$(x^{r}y)^{r}\in (x^{r+1},y^{r+1})(x^{r+1},y^{r+1},x^{r}y)^{r-1}$
for all $x,y\in A$.
\item[$(v)$] ${\rm dim}\, A\leq 1$.
\end{itemize}}

\noindent {\em Proof}. Implications $(i)\Rightarrow (ii)$ and
$(ii)\Rightarrow (iii)$ are obvious, $(iii)\Rightarrow (iv)$ is Remark
\ref{rest05} and $(iv)\Rightarrow (v)$ is Proposition \ref{rest1}. Let
us prove $(v)\Rightarrow (i)$. Let $A$ be an excellent ring of ${\rm
dim}\, A\leq 1$ and let $M$ be a finitely generated $A$-module.  If
${\rm dim}\, M=0$, then, by Proposition \ref{artim}, ${\rm
grt}(M)={\rm grt}_{A/{\rm Ann}_{A}(M)}(M)<\infty $. Therefore, we may
assume ${\rm dim}\, A=1$ and ${\rm dim}\, M=1$. Let ${\rm Min}(A)=\{
\mathfrak{p}_{1},\ldots ,\mathfrak{p}_{r}\} $ be the set of minimal
primes of $A$ and let ${\rm Ass}(A)={\rm Min}(A)\cup \{
\mathfrak{m}_{1},\ldots ,\mathfrak{m}_{s}\} $, $\mathfrak{m}_{i}\in
{\rm Max}(A)$, be the set of associated primes of $A$.  Since ${\rm
Ass}(A)$ is finite, by Propositions \ref{artim} and \ref{bo}, $\alpha
={\rm max}\{ {\rm grt}(M_{\mathfrak{p}})\mid \mathfrak{p}\in {\rm
Ass}(A)\} <\infty$.  Analogously, $\alpha ^{\prime}={\rm max}\{ {\rm
grt}(M_{\mathfrak{p}})\mid \mathfrak{p}\in {\rm Ass}(M)\} <\infty$.
If $r\geq 2$, and for each $1\leq i_{1}<\ldots <i_{l}\leq r$ with
$l\geq 2$, consider $\Gamma
_{i_{1},\ldots,i_{l}}=V(\mathfrak{p}_{i_{1}}+\ldots
+\mathfrak{p}_{i_{l}})$ and $\Gamma =\cup _{1\leq i_{1}<\ldots
<i_{l}\leq r}\Gamma _{i_{1},\ldots ,i_{l}}$. If $r=1$, take $\Gamma
=\emptyset$. In any case, $\Gamma $ is a closed finite subset of ${\rm
Spec}(A)$. By Proposition \ref{bo}, $\gamma ={\rm max}\{ {\rm
grt}(M_{\mathfrak{m}})\mid \mathfrak{m}\in \Gamma \} <\infty$.  Let
$\Sigma ={\rm Sing}(A/\mathfrak{p}_{1})\cup \ldots \cup {\rm
Sing}(A/\mathfrak{p}_{r})$, ${\rm Sing}(A/\mathfrak{p}_{i})=\{
\mathfrak{m}\in {\rm Max}(A)\mid \mathfrak{m}\supset \mathfrak{p}_{i}
\mbox{ and } A_{\mathfrak{m}}/\mathfrak{p}_{i}A_{\mathfrak{m}} \mbox{
is not regular }\} $. By hypothesis, ${\rm Sing}(A/\mathfrak{p}_{i})$
is a closed subset of ${\rm Spec}(A)$. In particular, ${\rm
Sing}(A/\mathfrak{p}_{i})$ and $\Sigma$ are finite. Again by
Proposition \ref{bo}, $\sigma ={\rm max}\{ {\rm
grt}(M_{\mathfrak{m}})\mid \mathfrak{m}\in \Sigma \} <\infty$. Now,
take $\mathfrak{m}\in {\rm Max}(A)$, $\mathfrak{m}\not\in {\rm
Ass}(A)\cup {\rm Ass}(M)\cup \Gamma \cup \Sigma$. Thus,
$A_{\mathfrak{m}}$ is a one dimensional Cohen-Macaulay local ring,
$M_{\mathfrak{m}}$ is a maximal Cohen-Macaulay
$A_{\mathfrak{m}}$-module, $\mathfrak{m}$ contains exactly one minimal
prime $\mathfrak{p}\in {\rm Min}(A)$, $\mathfrak{m}\supset
\mathfrak{p}$, and $A_{\mathfrak{m}}/\mathfrak{p}A_{\mathfrak{m}}$ is
a discrete valuation ring. Since $A$ is noetherian, there exists an
integer $n\geq 1$ such that $\mathfrak{n}(A)^{n}=0$. Thus
$\mathfrak{p}^{n}A_{\mathfrak{m}}=0$. By Lemma \ref{quiqm}, ${\rm
grt}(M_{\mathfrak{m}})\leq {\rm max}\{ n,\sum _{i=0}^{n-1}\mu
(\mathfrak{p}^{i}A_{\mathfrak{m}})\} \leq {\rm max}\{ n,\sum
_{i=0}^{n-1}\mu (\mathfrak{p}^{i})\}$. If $\mu =\sum _{i=0}^{n-1}\mu
(\mathfrak{p}^{i})$, then ${\rm grt}(M) ={\rm sup}\{ {\rm
grt}(M_{\mathfrak{p}})\mid \mathfrak{p}\in {\rm Spec}(A)\} \leq {\rm
max}\{ n,\mu ,\alpha ,\alpha ^{\prime}, \gamma ,\sigma \} <\infty
$. \fidemo

\begin{observacio} {\rm There exists (necessarily non $J-2$)
noetherian rings with ${\rm dim}\, A\leq 1$, but ${\rm grt}(A)=\infty
$.  For example, take $k$ a field and
$R=k[t_{1}^{2},t_{1}^{3},t_{2}^{3},t_{2}^{4},t_{2}^{5},\ldots ,
t_{g}^{g+1},t_{g}^{g+2},\ldots ,t_{g}^{2g+1},\ldots ]$. The ideals
$\mathfrak{p}_{g}=(t_{g}^{g+1},t_{g}^{g+2},\ldots ,t_{g}^{2g+1})$ are
prime of height 1. Let $S$ be the multiplicative closed set $S=R-\cup
\mathfrak{p}_{g}$ and $A=S^{-1}R$. Let
$\mathfrak{m}_{g}=S^{-1}\mathfrak{p}_{g}$. Since all prime ideals of
$R$ contained in $\cup \mathfrak{p}_{g}$ are contained in some
$\mathfrak{p}_{g}$, then $A$ is a one dimensional noetherian domain
with maximal ideals $\mathfrak{m}_{g}$ \cite{sv}. By Example
\ref{localtg}, ${\rm grt}(A_{\mathfrak{m}_{g}})=g+1$. Thus ${\rm
grt}(A)=\infty $. Remark ${\rm Sing}(A)={\rm Spec}(A)-\{ (0)\}$, so
$A$ is not $J-2$.}\end{observacio}

\noindent {\bf Theorem 1} {\em 
Let $A$ be an excellent (in fact $J-2$) ring and let $N\subseteq M$ be
two finitely generated $A$-modules such that ${\rm dim}(M/N)\leq 1$.
Then there exists an integer $s\geq 1$ such that, for all integers
$n\geq s$ and for all ideals $I$ of $A$,
\begin{eqnarray*}
I^{n}M\cap N=I^{n-s}(I^{s}M\cap N)\, .
\end{eqnarray*}}

\noindent {\em Proof}. Since ${\rm grt}_{A}(M/N)={\rm grt}_{A/J}(M/N)$
for $J={\rm Ann}_{A}(M/N)$, we can suppose that $A$ is an excellent
ring of ${\rm dim}\, A\leq 1$. Thus, by Theorems 2 and 3, $s(N,M)\leq
{\rm grt}(M/N)<\infty$. \fidemo

\bigskip

\noindent {\sc Acknowledgement}. I would like to thank J. \`Alvarez
and J. Masdemont for valuable comments. I wish to express my gratitude
to J.M. Giral and L. O'Carroll for their attention and interesting
conversations regarding this paper. This work was partially supported
by the UPC-PR9712 grant.

\end{document}